\documentclass{article}
\usepackage{graphicx}
\usepackage{amsmath,amssymb,latexsym}

\def\F{\mathcal{F}}
\def\tr{\textrm}

\def\dist{h}
\def\prf{\emph{Proof: }}
\def\qed{$\hfill\Box$}
\def\pic{\lambda}
\def\supp{\tr{supp}}
\newtheorem{theorem}{Theorem}
\newtheorem{thm}[theorem]{Theorem}
\newtheorem{question}[theorem]{Question}

\newtheorem{lem}[theorem]{Lemma}
\begin{document}
\title{Vertex Tur\'an problems in the hypercube}
\author{J Robert Johnson\thanks{School of Mathematical Sciences, Queen Mary University of London, E1 4NS, UK}\and John Talbot\thanks{
Department of Mathematics, University College London, WC1E 6BT, UK. 
Email: talbot@math.ucl.ac.uk.  This author is a Royal Society University Research Fellow.}}
\date\today
\maketitle
\begin{abstract} 
Let $\mathcal{Q}_n$ be the $n$-dimensional hypercube: the graph with vertex set $\{0,1\}^n$ and edges between vertices that differ in exactly one coordinate. For $1\leq d\leq n$ and $F\subseteq \{0,1\}^d$ we say that $S\subseteq \{0,1\}^n$ is  \emph{$F$-free} if every embedding $i:\{0,1\}^d\to \{0,1\}^n$ satisfies $i(F)\not\subseteq S$. We consider the question of how large $S\subseteq \{0,1\}^n$ can be if it is $F$-free.  In particular we generalise the main prior result in this area, for $F=\{0,1\}^2$, due to E.A.~Kostochka and prove a local stability result for the structure of near-extremal sets.

We also show that the density required to guarantee an embedded copy of at least one of a family of forbidden configurations may be significantly lower than that required to ensure an embedded copy of any individual member of the family.

Finally we show that any subset of the $n$-dimensional hypercube of positive density  will contain exponentially many points from some embedded $d$-dimensional subcube if $n$ is sufficiently large. 
\end{abstract}
\section{Introduction}
 For $n\geq 1$ let $V_n=\{0,1\}^n$. The \emph{$n$-dimensional hypercube}, $\mathcal{Q}_n$, is the graph with vertex set $V_n$ and  edges between vertices that differ in exactly one coordinate. 

An \emph{embedding} of $\mathcal{Q}_d$ into $\mathcal{Q}_n$ is an injective map $i:V_d\to V_n$ that preserves the edges of $\mathcal{Q}_d$. (Note that the image of $V_d$ under any such embedding consists of $2^{d}$ elements of $V_n$ given by fixing $n-d$ coordinates and allowing the other $d$ coordinates to vary.)

Given $F\subseteq V_d$, where $1\leq d\leq n$, we say that $S\subseteq V_n$ is  \emph{$F$-free} if every embedding $i:V_d\to V_n$ satisfies $i(F)\not\subseteq S$. 
For a family $\F$ of subsets of $V_d$ we say that $S$ is \emph{$\F$-free} if $S$ is $F$-free for all $F\in \F$. We define \[
\tr{exc}(n,\F)=\max\{|S|:S\subseteq V_n \tr{ is $\F$-free}\}.\]

It is easy to see (via averaging) that for any family $\F$ of subsets of $V_d$ the ratio $\tr{exc}(n,\F)/2^n$ is non-increasing and bounded below (by zero). Hence we can define the \emph{vertex Tur\'an density} by
\[
\pic(\F)=\lim_{n\to \infty}\frac{\tr{exc}(n,\F)}{2^n}.\]
We write $\pic(F)$ instead of $\pic(\{F\})$.

If $x\in V_n$ and $i\in [n]=\{1,2,\ldots,n\}$ then $x_i\in \{0,1\}$ denotes the $i$th coordinate of $x$. The \emph{support} of $x\in V_n$ is $\supp(x)=\{i\in [n]:x_i=1\}$ and the \emph{weight} of $x$ is $|x|=|\supp(x)|$. The set of all $x\in V_n$ of weight $r$ is called the $r$th \emph{layer}.

The quantity $\tr{exc}(n,V_2)$ was determined by Kostochka \cite{K}
(and independently by Johnson and Entringer \cite{EJ}). They also
showed that the unique largest $V_2$-free subset of $V_n$ can be obtained by deleting every third layer of $V_n$.  
\begin{thm}[Kostochka \cite{K}]\label{K:thm}
For $n\geq 2$ we have $\tr{exc}(n,V_2)=\lceil 2^{n+1}/3\rceil$. If $S\subseteq V_n$ is $V_2$-free and $|S|=\tr{exc}(n,V_2)$ then, up to automorphisms of $\mathcal{Q}_n$, $S$ is 
\[
S_{i}=\{x\in V_n: |x|\not\equiv i\ \tr{mod}\ 3\}\] for some $i\in \{0,1,2\}$.
\end{thm}
 
 The problem of determining $\tr{exc}(n,V_d)$ has been considered by various authors but mainly
when $d$ is close to $n$ (see \cite{KS}, \cite{R}, \cite{SB})  or when $n$ is small \cite{HH}.

Recently, Alon, Krech and Szab\'o \cite{AKS} described the problem of finding
$\pic(V_d)$. Their motivation was a related question of Erd\H os \cite{EC}  who asked for the
  largest number of edges in a $\mathcal{Q}_2$-free subgraph of
  $\mathcal{Q}_n$. The conjectured answer to this is $(1/2+o(1))e(\mathcal{Q}_n)$ while
  the best upper bound is around $0.62256 e(\mathcal{Q}_n)$ due to Thomason
  and Wagner \cite{TW} extending earlier work of Chung
  \cite{C}. A result relating the maximum density of edges of a
  $\mathcal{Q}_d$-free subgraph of $\mathcal{Q}_n$ to the analagous density for certain other
  forbidden subgraphs was proved by Offner \cite{OFF} using the supersaturation method. He also proved a vertex version of this result
  although our notion of containment as an embedded copy is slightly
  different from his. The general vertex version of the problem which we consider here is extremely
  natural but does not seem to have received attention.

\section{Results}
Our first result is a generalisation of Theorem \ref{K:thm}. We show that asymptotically the density required to guarantee a copy of $V_2$ is sufficient to ensure copies of other larger configurations. 

The configuration $G_d$ will be the set of all vertices of $V_d$ of weight zero or one together with a set of vertices of weight two whose supports form the edge set of a complete bipartite graph $K(\lceil d/2\rceil,\lfloor d/2\rfloor)$ (since all such configurations are isomorphic the precise choice of bipartition is unimportant, we will take the one given by parity of coordinates). Formally for $d \geq 2$ we define
\[
G_d=\{x\in V_d:|x|=0,1\tr{ or }(|x|=2, \supp(x)=\{i,j\},\ i\not\equiv j\tr{ mod $2$})\}.\] 
For example \[
G_2=V_2,\quad G_3=\{(0,0,0),(1,0,0),(0,1,0),(0,0,1),(1,1,0),(0,1,1)\}.\]
The relationship of our result to Theorem \ref{K:thm} is rather like that of the Erd\H os--Stone theorem for 3-chromatic graphs to Mantel's theorem in extremal graph theory.
\begin{figure}
\begin{center}
\includegraphics[scale=0.6]{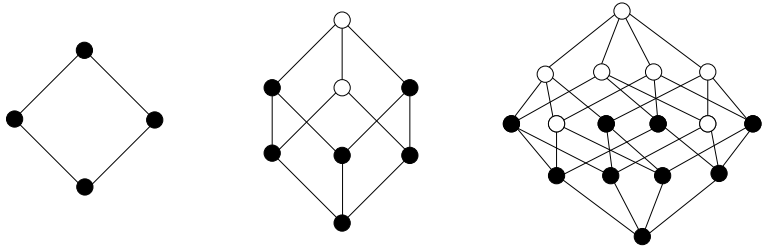}
\end{center}
\caption{Forbidden configurations $G_2$, $G_3$ and $G_4$ (black points)}
\end{figure}

 We also show that if $S\subseteq V_n$ is $G_d$-free and near-extremal in size then $S$ must locally resemble the ``two-out-of three-layers'' extremal construction given in Theorem \ref{K:thm}. 

We require some notation. Denote the \emph{Hamming distance} in $\mathcal{Q}_n$ by $\tr{dist}(x,y)$, this is the number of coordinates in which $x,y\in V_n$ differ.  For $l\geq 0$ and $x\in V_n$ let $\Gamma_l(x)=\{y\in V_n:\tr{dist}(x,y)=l\}$. Given $S\subseteq V_n$ we let $\dist_l(x)=|S\cap \Gamma_l(x)|$ denote the number of elements of $S$ at distance $l$ from $x$. 

Although $\dist_l(x)$ depends on the set $S$, for ease of notation we will suppress this. It will always be clear from the context what subset of $V_n$ we are considering.

 For $x\in V_n$, $S\subseteq V_n$ and $l\geq 1$ it is natural to view  $S\cap \Gamma_l(x)$ as an $l$-uniform hypergraph. Formally we define
\[
S^l(x)=\{\supp(x)\Delta\supp(y):y\in S\cap \Gamma_l(x)\}.\]
For a set $X$ and integer $r\geq 0$ we write $\binom{X}{r}$ for the family of all subsets of $X$ of size $r$.
So $A\in\binom{[n]}{l}$ belongs to $S^l(x)$ if flipping all of the coordinates of $x$ indexed by $A$  yields an element $y\in S$.

The precise definition of local stability is given in Theorem \ref{gd:thm} below but the three conditions may be paraphrased as follows.
\begin{itemize}
\item[(a)] For most $x\in V_n\setminus S$, most of the neighbours of $x$ (in $\mathcal{Q}_n$) belong to $S$.
\item[(b)] For most $x\in S$, approximately half of the neighbours of $x$ belong to $S$.
\item[(c)] For most $x\in S$ the graph $S^{(2)}(x)$ (corresponding to points at distance two from $x$ in $S$) is almost a clique on $[n]$ with a clique on $S^{(1)}(x)$ removed.\end{itemize}

\begin{thm}\label{gd:thm}
If $d\geq 2$ and $G_d$ is as defined above then
\begin{itemize}
\item[(i)] (Vertex Tur\'an density) \[\pic(G_d)=\frac{2}{3}.\]
\item[(ii)] (Local stability)  If $\epsilon>0$ there exists $\delta=\delta(\epsilon,d)$ satisfying $\lim_{\epsilon \to 0^+}\delta(\epsilon,d)=0$ and $n_0=n_0(\epsilon,d)$ such that if $n\geq n_0$ and $S\subseteq V_n$ is $G_d$-free with $|S|\geq (2/3-\epsilon)2^n$ then locally $S$ resembles the set $S_{0}=\{x\in V_n:|x|\not\equiv 0\ \tr{mod}\ 3\}$ in the following sense. There exists $T\subseteq  V_n$ with $|T|\leq \delta 2^n$ and
\begin{itemize} 
\item[(a)] $h_1(x)\geq (1-\delta)n$ for all $x\in V_n\setminus(S\cup T)$.
\item[(b)] $|h_1(x)-n/2|\leq \delta n$ for all $x\in S\setminus T$.
\item[(c)] $\left|S^2(x)\Delta\left(\binom{[n]}{2}\setminus \binom{S^1(x)}{2}\right)\right|\leq \delta \binom{n}{2} $ for all $x\in S\setminus T$.\end{itemize}
\end{itemize}
\end{thm}

Note that a ``global'' stability result cannot hold for this problem in the sense that there exist near-extremal size $G_d$-free subsets of $V_n$ that cannot be obtained from the ``two-out-of-three-layers'' construction by deleting/adding a small number of points and taking an automorphism of the hypercube. For example
\begin{multline*}
S=\{x\in V_n:|x|\leq n/2, |x|\not\equiv 0 \tr{ mod }3\}\ \cup\\\{x\in V_n:|x|\geq n/2+3, |x\Delta [n/2]|\not\equiv 0\tr{ mod }3\}.\end{multline*}

\medskip
Our second result (Theorem \ref{f12:thm}) shows that the density required to ensure a copy of at least one of a family of forbidden configurations may be significantly lower than that required to ensure a copy of any individual member of the family. This is in contrast to ordinary graph Tur\'an densities where the Erd\H os--Stone--Simonovits theorem implies that for any family $\F$ of graphs $\pi(\F)=\min\{\pi(F):F\in \F\}$ (where $\pi(\F)$ is the classical Tur\'an density of the family of graphs $\F$). This ``non-principality'' of the vertex Tur\'an density is analogous to that previously observed for $r$-uniform hypergraph Tur\'an densities by Balogh \cite{B} and Mubayi and Pikhurko \cite{MP} when $r\geq 3$.

\medskip
Let \[
F_1=\{(0,0,0),(1,0,0),(0,1,0),(0,0,1)\},\]\[F_2=\{(0,0,0),(1,1,0),(1,0,1),(0,1,1)\},\quad F_3=\{(0,0,0),(1,1,1)\}.\]
\begin{figure}
\begin{center}
\includegraphics[scale=0.35]{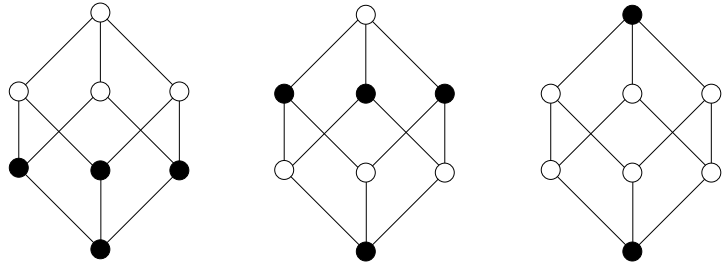}
\end{center}
\caption{Forbidden configurations $F_1$, $F_2$ and $F_3$}
\end{figure}
Theorem \ref{f12:thm} determines $\pic(\F)$ for all families $\F\subseteq \{F_1,F_2,F_3\}$.  In particular $\pic(\{F_2,F_3\})<\pic(\{F_1,F_2\})<\min\{\pic(F_1),\pic(F_2),\pic(F_3)\}$. 
\begin{thm}\label{f12:thm}
If $F_1$, $F_2$ and $F_3$ are as defined above then
\[
\pic(F_1)=\pic(F_2)=\pic(F_3)=\pic(\{F_1,F_3\})=\frac{1}{2},\]\[\pic(\{F_1,F_2\})=\frac{1}{3},\quad \pic(\{F_2,F_3\})=\pic(\{F_1,F_2,F_3\})=\frac{1}{4}.\]
\end{thm}

\medskip

Finally we consider the following question: given a subset of $V_n$ of positive density how many vertices must it contain from some $d$-dimensional subcube? For $1\leq t \leq d$ let $\F_{d,t}=\{F\subseteq V_d:|F|=t\}$. 
If $\lambda(\F_{d,t})=0$ then for any $\epsilon>0$ and $n\geq n_{0}(\epsilon,d)$ sufficiently large, any subset of $V_n$ of density at least $\epsilon$ will contain at least $t$ vertices from some $d$-dimensional subcube. We would like to determine the value of
\[
\mu(d)=\max\{t:\lambda(\F_{d,t})=0\}.\]
The following construction does not contain vertices from more than one layer of any $d$-dimensional subcube of $\mathcal{Q}_n$ and has density approximately $1/(d+1)$
\[
S_{d+1}=\{x\in V_n:|x|\equiv 0 \tr{ mod }d+1\}.\] Since the largest layer of $V_d$ has size $\binom{d}{\lfloor d/2\rfloor}$ this proves the upper bound in Theorem \ref{zd:thm}. The lower bound tells us that $\mu(d)$ is exponential in $d$.
\begin{thm}\label{zd:thm}
If $d\geq 2$ then \[
t_2(d)+t_3(d)\leq \mu(d)\leq \binom{d}{\lfloor d/2\rfloor},\] where
\[
t_2(d)=\left\{\begin{array}{ll}0,&\tr{if $\lceil d/3\rceil$ is odd},\\1,& \tr{otherwise}. \end{array}\right. \qquad t_3(d)=\left\{\begin{array}{ll}3^{d/3},& d\equiv 0 \mod 3,\\ 4\cdot 3^{(d-4)/3},& d\equiv 1 \mod 3,\\2\cdot 3^{(d-2)/3},& d\equiv 2 \mod 3.\end{array}\right.\]

\end{thm}

\section{Proofs}
We will make use of a number of classical results from extremal graph
and hypergraph theory. 
\begin{thm}[Mantel \cite{M}]\label{mantel:thm}
If $G=(V,E)$ is a triangle-free graph with $|V|=n$ then $|E|\leq n^2/4$.
\end{thm}

For $s\geq t\geq 1$ let $K(s,t)$ denote the complete bipartite graph with vertex classes of size $s$ and $t$.
For $r\geq 3$ and $t_1\geq t_2\geq\cdots\geq t_r\geq 1$ let $K^{(r)}(t_1,\ldots,t_r)$ denote the complete $r$-partite $r$-graph with vertex classes of size $t_1,\ldots,t_r$.
\begin{thm}[Erd\H os--Stone \cite{ESt}]\label{es:thm}
If $G=(V,E)$ is a $K(s,t)$-free graph and $|V|=n$ then $|E|=O(n^{2-1/t})$.
\end{thm}

\begin{thm}[Erd\H os \cite{E}]\label{erd:thm}
If the $r$-graph $G=(V,E)$ is $K^{(r)}(t_1,t_2,\ldots,t_r)$-free and $|V|=n$ then $|E|=O(n^{r-1/t_1})$.
\end{thm}

We will make repeated use of the following special case of the Cauchy--Schwarz inequality.
\begin{lem}\label{CS:lem} If $a_1,\ldots,a_s\in \mathbb{R}$ and $\frac{1}{s}\sum_{i=1}^s a_i\geq A$ then
\[
\frac{1}{s}\sum_{i=1}^s a_i^2\geq A^2.\]\end{lem}


 For $x,y\in V_n$ let $\dist_l(x,y)$ denote the number of elements of $S$ at distance $l$ from both $x$ and $y$, i.e.
\[
\dist_l(x,y)=|S\cap\Gamma_l(x)\cap\Gamma_l(y)|.\]
%
%

The following simple lemma underpins all our results.
\begin{lem}\label{key:lem}
If $S\subseteq V_n$ and $l\geq 1$ then
\begin{itemize}
\item[(i)]
\[
\sum_{v\in V_n}\dist_l^2(v)= \sum_{x\in S}\binom{2l}{l}\dist_{2l}(x)+O(n^{2l-1}2^n).
\]
\item[(ii)]
\[
\sum_{v\in S}\dist_l^2(v)= \sum_{x\in S}\sum_{z\in S\cap\Gamma_{2l}(x)}\dist_l(x,z)+O(n^{2l-1}2^n).\]
\end{itemize}
\end{lem}
\prf
 For (i) consider the sum
\[
\sum_{x\in S}\sum_{y\in \Gamma_l(x)}\dist_l(y).\]
For each $v\in V_n$ the term $\dist_l(v)$ occurs once for each element of $S$ at distance $l$ from $v$, i.e.~$\dist_l(v)$ times. Hence 
\[
\sum_{x\in S}\sum_{y\in \Gamma_l(x)}\dist_l(y)=\sum_{v\in V_n}\dist_l^2(v).\]
Moreover for a fixed choice of $x\in S$ the inner sum counts elements of $S$  that can be reached from $x$ by flipping $l$ coordinates of $x$ and then flipping $l$ coordinates of the resulting point. Hence if $0\leq k\leq l$ then the sum counts those elements of $S$ at distance $2k$ from $x$ precisely $\binom{2k}{k}\binom{n-2k}{l-k}$ times. Thus
\begin{eqnarray*}
\sum_{x\in S}\sum_{y\in \Gamma_l(x)}\dist_l(y)&=&\sum_{x\in S}\sum_{k=0}^l\binom{2k}{k}\binom{n-2k}{l-k}\dist_{2k}(x) \\
 & = & \sum_{x\in S}\binom{2l}{l}\dist_{2l}(x)+O(n^{2l-1}2^n),\end{eqnarray*}
since $\dist_{2k}(x)=O(n^{2k})$ for any $0\leq k \leq l$ and $|S|\leq 2^n$. Hence (i) holds.

For  (ii) consider the sum
\[
\sum_{x\in S}\sum_{y\in S\cap \Gamma_l(x)}\dist_l(y).\]
For each $v\in S$ the term $\dist_l(v)$ occurs once for each element of $S$ at distance $l$ from $v$, i.e.~$\dist_l(v)$ times. Hence 
\[
\sum_{x\in S}\sum_{y\in S\cap\Gamma_l(x)}\dist_l(y)=\sum_{v\in S}\dist_l^2(v).\]
The same argument as used for (i) implies that the contribution to the LHS of this sum from $z\in S$ satisfying $\tr{dist}(x,z)<2l$ is at most $O(n^{2l-1}2^n)$. Finally $z\in S\cap \Gamma_{2l}(z)$ contributes one to this sum for each choice of $y\in S\cap\Gamma_l(x)$ such that $\tr{dist}(y,z)=l$, i.e.~$\dist_l(x,z)$ times. The result follows.
\qed

\bigskip
\emph{Proof of Theorem \ref{f12:thm}: }
For the lower bounds note that the following sets are $F_1$, $F_2$, $\{F_1,F_2\}$ and $\{F_2,F_3\}$-free respectively:
\[
S_1=\{x\in V_n:|x|\equiv 0\tr{ mod } 2\},\quad S_2=\{x\in V_n:|x|\equiv 0,1 \tr{ mod } 4\},\] \[S_{1,2}=\{x\in V_n:|x|\equiv 0 \tr{ mod } 3\}\quad S_{2,3}=\{x\in V_n:|x|\equiv 0 \tr{ mod } 4\}.\]
 Moreover these sets have asymptotic densities $1/2$, $1/2$, $1/3$ and $1/4$ respectively. Since $S_1$ is also $F_3$-free and $S_{2,3}$ is also $F_1$-free it is sufficient to prove that these values are also upper bounds for the vertex Tur\'an densities.

Let $T_1$ be $F_1$-free with $|T_1|=\alpha_12^n$. If $x\in T_1$ then $\dist_1(x)=|T_1\cap\Gamma_1(x)|\leq 2$, hence 
\begin{eqnarray*}
n\alpha_12^n&=&\sum_{x\in V_n}\dist_1(x)\\
& \leq & \sum_{x\in T_1}2+\sum_{x\in V_n\setminus T_1}n\\
& = &2\alpha_1 2^n+n(1-\alpha_1)2^n.\end{eqnarray*}
So $\alpha_1\leq n/(2n-2)$ and hence $\pic(F_1)=1/2$. 

Similarly if $T_3$ is $F_3$-free and $|T_3|=\alpha_32^n$ then $h_3(x)=|T_3\cap \Gamma_3(x)|=0$ for all $x\in T_3$. Hence
\[
\binom{n}{3}\alpha_3 2^n=\sum_{x\in V_n}\dist_3(x)\leq \binom{n}{3}|V_n\setminus T_3|=
\binom{n}{3}(1-\alpha_3)2^n.\]
Thus $\alpha_3\leq 1/2$ and hence $\pic(F_3)= 1/2$.

Let $T_2$ be $F_2$-free with $|T_2|=\alpha_2 2^n$. For $x\in T_2$ let \[
T_2^2(x)=\{\supp(x)\Delta\supp(y):y\in T_2\cap \Gamma_2(x)\}.\]
Consider the graph with vertex set $[n]$ and edge set $T_2^2(x)$. Since $T_2$ is $F_2$-free this graph is triangle-free (a triangle would correspond to $a,b,c\in T_2$ such that $\{x,a,b,c\}$ forms a copy of $F_2$ in $T_2$). Hence, by Mantel's theorem, $|T^2_2(x)|\leq n^2/4$. Thus for $x\in T_2$ we have $\dist_2(x)=|T_2^2(x)|\leq n^2/4$, so Lemma \ref{key:lem} (i) with $l=1$ implies that \begin{equation}\label{1:eq}
\sum_{x\in V_n}\dist_1^2(x)\leq n^2\alpha_2 2^{n-1}+O(n2^n).\end{equation}
Since $\sum_{x\in V_n}\dist_1(x)=n\alpha_22^{n}$, Cauchy--Schwarz implies that\[ n^2\alpha_2^2 2^n\leq n^2 \alpha_22^{n-1}+O(n2^n).\]
Hence $\alpha_2\leq 1/2+o(1)$ and so $\pic(F_2)=1/2$.

Now let $T_{1,2}$ be $\{F_1,F_2\}$-free with $|T_{1,2}|=\alpha_{1,2}2^n$. For each $x\in T_{1,2}$ we have both $\dist_1(x)\leq 2$ and $\dist_2(x)\leq n^2/4$. So (\ref{1:eq}) holds with $\alpha_2$ replaced by $\alpha_{1,2}$ and 
\[\sum_{x\in V_n \setminus T_{1,2}}\dist_1(x) \geq (n-2)\alpha_{1,2}2^n.\]
Hence Cauchy--Schwarz implies that
\[
\left(\frac{(n-2)\alpha_{1,2}}{1-\alpha_{1,2}}\right)^2(1-\alpha_{1,2})2^n\leq n^2\alpha_{1,2}2^{n-1}+O(n2^n).\]
So $\alpha_{1,2}\leq n^2/(2(n-2)^2+n^2)+o(1)$ and $\pic(\{F_1,F_2\})= 1/3$. 

Finally let $T_{2,3}$ be $\{F_2,F_3\}$-free and $|T_{2,3}|=\alpha_{2,3}2^n$. Since $T_{2,3}$ is $F_2$-free, (\ref{1:eq}) holds with $\alpha_2$ replaced by $\alpha_{2,3}$.

Let $Y=\{x\in V_n:\dist_1(x)\geq 3\}$ and $|Y|=\beta 2^n$.
Since $\sum_{x\in V_n}h_1(x)=n\alpha_{2,3}2^n$ and $h_1(x)\leq 2$ for $x\in V_n \setminus Y$ we have
\[
\sum_{x\in Y}h_1(x)\geq n\alpha_{2,3}2^n -2(1-\beta)2^n.\] Hence, using Cauchy--Schwarz,
\[
\sum_{x\in V_n}h_1^2(x)\geq \sum_{x\in Y}h_1^2(x)\geq \frac{(\alpha_{2,3}n-2+2\beta)^2 2^n}{\beta}.\]
Thus
\[
\alpha_{2,3}\leq \frac{\beta}{2}+o(1).\]
If we show that $\beta \leq 1/2$ we will be done. If $\beta>1/2$ then $|Y|>2^{n-1}$ and so there exist $y,z\in Y$ such that $\tr{dist}(y,z)=1$. Let $a,b,c\in T_{2,3}\setminus \{y,z\}$ satisfy $\tr{dist}(a,y)=\tr{dist}(b,z)=\tr{dist}(c,z)=1$ (such points exist by the definition of $Y$). Now either $\tr{dist}(a,b)=3$ or $\tr{dist}(a,c)=3$ and so $T_{2,3}$ contains a copy of $F_3$. Hence $\beta\leq 1/2$ and so $ \alpha_{2,3}\leq 1/4+o(1)$ and $\pic(\{F_2,F_3\})=1/4$.
\qed

\bigskip
For $d\geq 3$ define
\[
F_1^d=\{x\in V_d:|x|=0,1\},\quad F_d^d=\{x\in V_d:|x|=0,d\}.\]
\[
F_2^d=\{x\in V_d:|x|=0 \tr{ or }(|x|=2, \supp(x)=\{i,j\},\ i\not\equiv j\tr{ mod $2$})\}.\] 
The proof of Theorem \ref{f12:thm} is easily extended to give the following result.
\begin{thm}
If $d_1,d_2,d_3\geq 3$ with $d_3$ odd then
\[
\pic(F^{d_1}_1)=\pic(F^{d_2}_2)=\pic(F^{d_3}_{d_3})=\pic(\{F^{d_1}_1,F^{d_3}_{d_3}\})=\frac{1}{2},\]\[\pic(\{F^{d_1}_1,F^{d_2}_2\})=\frac{1}{3},\quad \pic(\{F^{d_2}_2,F^{d_3}_{d_3}\})=\pic(\{F^{d_1}_1,F^{d_2}_2,F^{d_3}_{d_3}\})=\frac{1}{4}.\]

\end{thm}
In fact it is an immediate consequence of a result of Chung, F\"uredi, Graham and Seymour \cite{Cet} that $\tr{exc}(n,F_1^d)=2^{n-1}$ for $n$ sufficiently large.

\bigskip
\emph{Proof of Theorem \ref{gd:thm}:}
We start by proving $\lambda(G_d)=2/3$. Since \[
S_0=\{x\in V_n:|x|\not \equiv 0 \tr{ mod }3\}\] is $G_d$-free we have $\lambda(G_d)\geq 2/3$.

Let $2\leq d\leq n$ and $S\subseteq V_n$ be $G_d$-free. If $|S|=\alpha2^n$ then we need to show that $\alpha\leq 2/3+o(1)$. 

For $x\in S$ and $l\geq 1$ recall that \[
S^l(x)=\{\supp(x)\Delta\supp(y):y\in S\cap \Gamma_l(x)\}.\]
The fact that $S$ is $G_d$-free implies that for any $x\in S$ the graph with vertex set $S^1(x)$ and edge set $S^2(x)$ is $K(\lceil d/2\rceil,\lfloor d/2\rfloor)$-free and hence $K(d,d)$-free. Thus the Erd\H os--Stone theorem implies that it contains at most $O(n^{2-1/d})$ edges. Hence
\begin{equation}\label{h2:eq}
\dist_2(x)\leq \binom{n}{2}-\binom{\dist_1(x)}{2}+O(n^{2-1/d}).\end{equation}
Applying Lemma \ref{key:lem} (i) with $l=1$ we obtain
\begin{equation}\label{gd:eq}
\sum_{x\in V_n}\dist_1^2(x)\leq\sum_{x\in S}( n^2-\dist_1^2(x)) + O(n^{2-1/d}2^n).\end{equation}
Now let $\beta$ be defined by $\sum_{x\in S}\dist_1(x)=\beta \alpha n2^n$, so $0\leq \beta \leq 1$. Using (\ref{gd:eq}) and applying Cauchy--Schwarz to the sums $\sum_{x\in S}h_1^2(x)$ and $\sum_{x\in V_n\setminus S}h_1^2(x)$ we obtain
\[
2\beta^2n^2 \alpha 2^{n}+\left(\frac{(1-\beta)\alpha n}{1-\alpha}\right)^2(1-\alpha)2^n\leq \alpha n^2 2^n +O(n^{2-1/d}2^n).\]
So \[
\alpha(2-2\beta -\beta^2)\leq 1-2\beta^2+o(1).\]
If $2-2\beta-\beta^2>0$ then \begin{equation}\label{alpha:eq}
\alpha\leq \frac{1-2\beta^2}{2-2\beta-\beta^2}+o(1)\end{equation} and the RHS is maximised at $\beta=1/2$ when it equals $2/3+o(1)$. So  suppose that $2-2\beta -\beta^2\leq 0$. Since $0\leq \beta \leq 1$ this implies that $\beta \geq \sqrt{3}-1$.

 If $x\in S$ and $z\in S\cap \Gamma_2(x)$ then $\dist_1(x,z)=|S\cap\Gamma_1(x)\cap\Gamma_1(z)|\leq 2$. Moreover since $S$ is $G_d$-free the Erd\H os--Stone theorem implies that \[|\{z\in S\cap \Gamma_2(x): \dist_1(x,z)=2\}|=O(n^{2-1/d}).\]
Finally \[|\{z\in S\cap \Gamma_2(x): \dist_1(x,z)=1\}|\leq \binom{n}{2}-\binom{\dist_1(x)}{2}.\]
Hence Lemma \ref{key:lem} (ii) with $l=1$ implies that
\[
\sum_{x\in S}\dist_1^2(x)\leq \sum_{x\in S}\left(\binom{n}{2}-\binom{\dist_1(x)}{2}\right)+O(n^{2-1/d}2^n).\]
Using $\sum_{x\in S}\dist_1(v)=\beta\alpha n2^n$ and Cauchy--Schwarz we obtain
\[
3\beta^2n^2\alpha 2^{n-1}\leq n^2\alpha 2^{n-1}+O(n^{2-1/d}2^n).\]
Hence $\beta \leq 1/\sqrt{3}+o(1)<\sqrt{3}-1$ for $n$ large, and so $\pic(G_d)=2/3$. 

\medskip
We now need to show that the local stability conditions hold. Suppose that $S\subseteq V_n$ is $G_d$-free and has size $|S|\geq (2/3-\epsilon)2^n$ for some $\epsilon>0$. For $n$ large (\ref{alpha:eq}) implies that $|S|\leq (2/3+\epsilon)2^n$. If $\epsilon\geq 1/100$ then we may take $\delta=1$ and the conditions hold trivially so suppose that $\epsilon\leq 1/100$.

Since $\alpha\geq 2/3-\epsilon$,   (\ref{alpha:eq}) implies that $|\beta-1/2|\leq\sqrt{\epsilon}$ for $n$ large. Let $\delta_1=2\epsilon^{1/4}$ and suppose there exists $W\subset V_n\setminus S$ such that $h_1(x)<(1-\delta_1)n$ for all $x \in W$ and $|W|\geq \delta_12^n$. Then \begin{eqnarray*}
\sum_{x\in V_n\setminus S}h_1(x)&\leq& |V_n\setminus(S\cup W)|n+(1-\delta_1)n|W|\\ &\leq &\left(1-\alpha-2\epsilon^{1/4}\right)n2^n+2\epsilon^{1/4}(1-2\epsilon^{1/4})n2^n\\
&\leq & \left(\frac{1}{3}-3\sqrt{\epsilon}\right)n2^n\end{eqnarray*}
But \begin{equation}\label{sc2:eq}
\sum_{x\in V_n\setminus S}h_1(x)=(1-\beta)\alpha n 2^n\geq \left(\frac{1}{3}-\frac{7\sqrt{\epsilon}}{6}\right)n2^n.\end{equation}
Hence (a) holds for any $\delta\geq 2\epsilon^{1/4}$.

For (b) we will require the following defect form of the Cauchy--Schwarz inequality (cf.~Bollob\'as \cite{Bbook} page 125).
\begin{lem}\label{cs:lem}
If $a_1,\ldots,a_s\in \mathbb{R}$, $1\leq t\leq s$, $\frac{1}{s}\sum_{i=1}^s a_i=A$, $\frac{1}{t}\sum_{i=1}^t a_i\geq A'$, $t\geq \gamma s$ and $A'\geq A+\eta$ then \[
\frac{1}{s}\sum_{i=1}^s a_i^2\geq A^2+\gamma \eta^2.\]
\end{lem}
Let $\delta_2=2\epsilon^{1/6}$ and suppose there exists $W\subset S$ satisfying $|W|\geq \delta_22^n$ and $|h_1(x)-n/2|\geq\delta_2 n$ for all $x\in W$. Let
\[
W^+=\{x\in W: h_1(x)\geq n/2+\delta_2n\}\] and $W^-=W\setminus W^+$. Suppose that $|W^+|\geq \delta_2 2^{n-1}$ (the case $|W^-|\geq \delta_22^{n-1}$ is similar so we omit it). 

Now \[
\frac{1}{|S|}\sum_{x\in S}h_1(x)=\beta n,\quad\frac{1}{|W^+|}\sum_{x\in W^+}h_1(x)\geq \frac{n}{2}+\delta_2n\]
so using $|\beta-1/2|<\sqrt{\epsilon}$ we have
\[
\frac{1}{|W^+|}\sum_{x\in W^+}h_1(x) \geq \frac{1}{|S|}\sum_{x\in S}h_1(x)+n(\delta_2-\sqrt{\epsilon}).\]
Since $|W^+|\geq \delta_2 2^{n-1} \geq \delta_2 |S|/2$, Lemma \ref{cs:lem} (with $A=\beta n$, $\eta=n(\delta_2-\sqrt{\epsilon})$ and $\gamma=\delta_2/2$) implies that 
\begin{equation}\label{s:eq}
\frac{1}{|S|}\sum_{x\in S}h_1^2(x) \geq \beta^2n^2 + \frac{\delta_2 n^2(\delta_2-\sqrt{\epsilon})^2}{2}.\end{equation}
Using (\ref{sc2:eq}) and Cauchy--Schwarz (Lemma \ref{CS:lem}) we also have
\begin{equation}\label{sc:eq}
\sum_{x\in V_n\setminus S}h_1^2(x)\geq \frac{1}{3}\left(1-7\sqrt{\epsilon}\right)n^22^n.\end{equation}
Now (\ref{gd:eq}) implies that for $n$ large
\[
\sum_{x\in S}2h_1^2(x)+\sum_{x\in V_n\setminus S}h_1^2(x)\leq n^2|S|+\epsilon n^22^n.\]
Using (\ref{s:eq}) and (\ref{sc:eq}) this yields
\[
\left(\frac{2}{3}-\epsilon\right)\left(2\beta^2+\delta_2(\delta_2-\sqrt{\epsilon})^2\right)+\frac{1}{3}\left(1-7\sqrt{\epsilon}\right)\leq \frac{2}{3}+2\epsilon.\] Substituting $\delta_2=2\epsilon^{1/6}$ and using $\beta\geq 1/2-\sqrt{\epsilon}$, $\epsilon\leq 1/100$ we obtain a contradiction. Hence (b) holds for any $\delta\geq 2 \epsilon^{1/6}$.

Finally for (c) let $\delta_3=4\epsilon^{1/4}$. Since $S$ is $G_d$-free the Erd\H os--Stone theorem implies that for any $x\in S$, $S^2(x)$ contains at most $O(n^{2-1/d})$ edges from $\binom{S^1(x)}{2}$. Hence
\[
\left|S^2(x)\cap \binom{S^1(x)}{2}\right|=O(n^{2-1/d})\leq \frac{\delta_3}{2}\binom{n}{2},\] for $n$ large. So suppose there exists $W\subseteq S$ such that $|W|\geq \delta_3 2^n$ and for all $x\in W$
\[
\left|\left(\binom{[n]}{2}\setminus\binom{S^1(x)}{2}\right)\setminus S^2(x)\right|\geq \frac{\delta_3}{2}\binom{n}{2}.\] 
Since (\ref{h2:eq}) holds for all $x\in S$, Lemma \ref{key:lem} (i) with $l=1$ implies that
\[
\sum_{x\in V_n}h_1^2(x)\leq \sum_{x\in S}(n^2-h_1^2(x)+O(n^{2-1/d}))-\delta^2_3 \binom{n}{2}2^n.\]
So for $n$ large
\[
2\sum_{x\in S}h_1^2(x)+\sum_{x\in V_n\setminus S}h_1^2(x)\leq n^2 |S|-\frac{\delta_3^2}{3}n^22^n.\] Applying Cauchy--Schwarz to $\sum_{x\in S}h_1^2(x)$ and using (\ref{sc:eq}) we obtain
\[
2\beta^2\alpha + \frac{1}{3}\left(1-7\sqrt{\epsilon}\right)\leq \alpha-\frac{\delta_3^2}{3}.\]
However this gives a contradiction using $|\beta-1/2|\leq \sqrt{\epsilon}$, $\alpha \leq 2/3+\epsilon$ and $\delta_3 =4\epsilon^{1/4}$. Hence (c) holds for any $\delta\geq 4\epsilon^{1/4}$. Thus we can satisfy all of the local stability conditions by taking $\delta=4\epsilon^{1/6}$, which clearly also satisfies the condition $\lim_{\epsilon\to 0^+}\delta(\epsilon)=0$.\qed

\bigskip
\emph{Proof of Theorem \ref{zd:thm}:}
Let $d\geq 2$ and $\epsilon>0$ be given. Let $n$ be large and $S\subseteq V_n$ satisfy $|S|\geq \epsilon 2^n$. For $r \geq 1$ we have
\[
\sum_{x\in V_n}\dist_r(x)=\binom{n}{r}|S|.\]
Hence $|S|\geq \epsilon 2^n$ implies that for any $r\geq 1$ there exists $x\in V_n$ such that $\dist_r(x)\geq \epsilon\binom{n}{r}$. Let $r=\lceil d/3\rceil$ and $3\geq p_1\geq p_2\geq\cdots\geq p_r\geq 2$ satisfy $\sum_{i=1}^rp_i=d$. If $n$ is sufficiently large relative to $\epsilon$ and $d$ then Theorem \ref{erd:thm} implies that the $r$-graph
\[
S^r(x)=\{\supp(x)\Delta \supp(y):y\in S\cap \Gamma_r(x)\}\]
contains a copy of $K:=K^{(r)}(p_1,p_2,\ldots,p_r)$. It is easy to check that this $r$-graph has $t_3$ edges (where $t_3$ is defined in the statement of Theorem \ref{zd:thm}). Moreover, since $K$ is an $r$-graph with $d$ vertices, the corresponding subset of $S$ lies in a copy of $V_d$. Thus $\lambda(\F_{d,t_3})=0$.

If $r=\lceil d/3\rceil$ is odd then $t_2=0$ and the proof is complete so suppose that $\lceil d/3\rceil$ is even. Now Lemma \ref{key:lem} (i) (with $l=r/2$) followed by an application of Cauchy--Schwarz implies that
\begin{eqnarray*}
\binom{r}{r/2}\sum_{x\in S}\dist_r(x)&=& \sum_{x\in V_n}h_{r/2}^2(x)+O(n^{r-1}2^n)\\
&\geq & \binom{n}{r/2}^2\epsilon^2 2^n+O(n^{r-1}2^n).\end{eqnarray*}
Hence there is $x\in S$ such that
\[
\dist_r(x)\geq \epsilon \binom{n}{r}+O(n^{r-1}).\] The same argument as above implies that $S^{r}(x)$ contains a copy of $K$. Hence there is a copy of $V_d$ containing $t_3+1$ points from $S$, so $\lambda(\F_{d,t_3+t_2})=0$.\qed 

\section{Questions}
There are many open problems concerning the vertex Tur\'an density. We collect what seem to be the most appealing ones here.

All of the constructions we have considered are of the form $\{x\in V_n : |x|\in I\}$ for some $I\subseteq[n]$.
\begin{question}\label{layers:q}
Is it true that for any  family $\F=\{F_1,\ldots,F_k\}$ of subsets of $V_d$ there are sets $I_n\subseteq[n]$ with
\[
S_{n}=\{x\in V_n : |x|\in I_n\}
\]
$\F$-free for all $n$ and 
\[
\lim_{n\to \infty}\frac{|S_{n}|}{2^n}=\pic(\F)?
\]
\end{question}

All of our results are for the vertex Tur\'an density, one could also ask for the exact value of $\tr{exc}(n,\F)$ and what the extremal examples are. 

As we mentioned above the question of determining (or improving bounds on) $\pic(V_d)$ was posed by Alon, Krech and Szab\'o \cite{AKS}. This is perhaps the most natural forbidden configuration to consider. 
%
%

Recall that $\F_{d,t}=\{F\subseteq V_d:|F|=t\}$ and $\mu(d)=\max\{t:\lambda(\F_{d,t})=0\}$. By Theorem \ref{zd:thm} we have $\mu(d)\leq \binom{d}{\lfloor d/2\rfloor}$. 
\begin{question}\label{zd:q}
Is it true that for all $d\geq 2$, $\mu(d)=\binom{d}{\lfloor d/2\rfloor}$?
\end{question}
Theorem \ref{zd:thm} tells us that $\mu(3)=3$ and $\mu(4)\geq
5$. Whether or not $\mu(4)=6$ is unresolved.

\section{Acknowledgement}
We thank the two anonymous referees for a number of comments, corrections and suggestions which greatly improved the presentation.

\end{document}